\documentclass[12pt]{amsart}
\usepackage[utf8]{inputenc}
\usepackage[english]{babel}
\usepackage{amsmath,amssymb,a4wide,latexsym,amscd,amsthm}
\usepackage{graphics}

\theoremstyle{plain}
\newtheorem{theorem}{Theorem}[section]
\newtheorem*{A}{Theorem A}
\newtheorem*{B}{Theorem B}
\newtheorem*{C}{Theorem C}
\newtheorem{lemma}{Lemma}[section]

\theoremstyle{definition}

\newtheorem{remark}{Remark}[section]

\numberwithin{equation}{section}

\begin{document}

\title[]
{When is the sum of closed subspaces of a Hilbert space closed?}

\author[]{Ivan Feshchenko}
\address{Institute of Mathematics of NAS of Ukraine, Kyiv, Ukraine}
\email{ivanmath007@gmail.com}

\begin{abstract}
We provide a sufficient condition for a finite number of closed subspaces
of a Hilbert space to be linearly independent
and their sum to be closed.
Under this condition a formula for the orthogonal projection
onto the sum is given.
We also show that this condition is sharp
(in a certain sense).
\end{abstract}

\subjclass[2010]{Primary 46C05, 46C07; Secondary 47B15}

\keywords{Hilbert space, closed subspace, sum of subspaces,
orthogonal projection.}

\maketitle

\section{Introduction}

\subsection{A few auxiliary notions}

\textbf{Complemented subspaces in Banach spaces.}
Let $X$ be a real or complex Banach space.
By a subspace of $X$ we will mean a linear subset of $X$.
Let $M$ be a subspace of $X$.
$M$ is said to be complemented in $X$ 
if there exists a subspace $N$ (a complement) such that 
$X$ is the topological direct sum of $M$ and $N$.
This means that the sum operator $S:M\times N\to X$ 
defined by $S(x,y)=x+y$, $x\in M, y\in N$ is an isomorphism 
(of normed linear spaces).
Here $M\times N$ is the linear space of all pairs $(x,y)$ 
with $x\in M, y\in N$ endowed with the norm $\|(x,y)\|=\|x\|+\|y\|$.
One can easily check that $M$ is complemented in $X$ if and only if 
there exists a continuous linear projection onto $M$, i.e.,
a continuous linear operator $P:X\to X$ such that 
$Px\in M$ for all $x\in X$ and $Px=x$ for $x\in M$.

Each complemented subspace is closed 
(this follows from the fact that $M=S(M\times\{0\})$).
Note that one can give the following (equivalent) 
definition of complementability:
a subspace $M$ is said to be complemented in $X$ if $M$ is closed 
and there exists a closed subspace $N$ such that 
$M\cap N=\{0\}$ and $M+N=X$
(the equivalence of the definition to the original 
follows from the fact that each complemented subspace is closed and
the Banach inverse mapping theorem).

If $X$ is a Hilbert space, 
then each closed subspace $M$ of $X$ is complemented in $X$
(one can consider the orthogonal decomposition $X=M\oplus M^\bot$ 
or, equivalently, the orthogonal projection onto $M$).

\textbf{Sum of subspaces.}
Let $V$ be a vector space and $V_1,...,V_n$ be subspaces of $V$.
Define the sum of $V_1,...,V_n$ in the natural way, namely,
$$
V_1+...+V_n:=\{x_1+...+x_n\,|\,x_1\in V_1,...,x_n\in V_n\}.
$$
It is clear that $V_1+...+V_n$ is a subspace of $V$.

\textbf{Linear independence.}
Let $V$ be a vector space and $V_1,...,V_n$ be subspaces of $V$.
The system of subspaces $V_1,...,V_n$ is said to be linearly independent 
if an equality $x_1+...+x_n=0$, where $x_1\in V_1,...,x_n\in V_n$,
implies that $x_1=...=x_n=0$.

\subsection{Notation}
Throughout the paper, $X$ is a real or complex Banach space 
with norm $\|\cdot\|$.
When $X$ is a Hilbert space we denote by $\langle\cdot,\cdot\rangle$
the inner product in $X$.
The identity operator on $X$ is denoted by $I$ 
(throughout the paper it is clear which Banach space is being considered).
All operators in the paper are continuous linear operators.
In particular, by a projection we always mean a continuous linear projection.
The kernel and range of an operator $T$ 
will be denoted by $\ker(T)$ and $Ran(T)$, respectively.
For a continuous linear operator $T$ between two Hilbert spaces
we denote by $T^*$ its adjoint.
All vectors are vector-columns; 
the superscript "t" means transpose.

\subsection{}
Starting point for this paper are the main results of our paper \cite{F20}
where the following questions are studied.
Let $X$ be a Banach space and $X_1,...,X_n$ be complemented subspaces of $X$.

\textbf{Question 1:}
\textit{Is $X_1+...+X_n$ complemented in $X$?}

If Question 1 has positive answer (for given $X_1,...,X_n$),
then the next natural question arises:

\textbf{Question 2:}
\textit{Suppose that we know some projections 
$P_1,...,P_n$ onto $X_1,...,X_n$, respectively.
Is there a formula for a projection onto $X_1+...+X_n$ 
(in terms of $P_1,...,P_n$) (of course, under certain conditions)?}

\subsection{}
Let us present the main results of \cite{F20}.
First, we provided a sufficient condition for the sum
of complemented subspaces to be complemented.

\begin{A}(\cite[Theorem 2.1]{F20})
Let $X$ be a Banach space, 
$X_1,...,X_n$ be complemented subspaces of $X$ and 
$P_1,...,P_n$ be projections onto $X_1,...,X_n$, respectively.
Let $\varepsilon_{ij}$, $i\neq j$, $i,j\in\{1,...,n\}$,
be nonnegative numbers such that
\begin{equation*}
\|P_i x\|\leqslant\varepsilon_{ij}\|x\|,\quad x\in X_j,
\end{equation*}
for any distinct $i,j\in\{1,...,n\}$.
Define the $n\times n$ matrix $E=(e_{ij})$ by
$$e_{ij}=
\begin{cases}
0, &\text{if $i=j$,}\\
\varepsilon_{ij}, &\text{if $i\neq j$.}
\end{cases}
$$
Denote by $r(E)$ the spectral radius of $E$ and set $A:=P_1+...+P_n$.

If $r(E)<1$, 
then the subspaces $X_1,...,X_n$ are linearly independent, 
their sum is complemented in $X$, and 
$\ker(P_1)\cap...\cap\ker(P_n)$ is a complement for $X_1+...+X_n$ in $X$.
Moreover, $I-(I-A)^N$ converges uniformly to the projection $P$ 
onto $X_1+...+X_n$ along $\ker(P_1)\cap...\cap\ker(P_n)$ as $N\to\infty$.
\end{A}

\subsection{}
The next result shows that the rate of convergence of
$I-(I-A)^N$ to $P$ can be estimated from above by $C\alpha^N$, 
where $\alpha\in[0,1)$.
To formulate the result we need the following notation: 
for two vectors $u,v\in\mathbb{R}^n$ we will write $u\leqslant v$ 
if $u\leqslant v$ coordinatewise.
To make the formulation of the result clearer 
we first make the following important remark.
Since $E$ is a nonnegative matrix, the condition $r(E)<1$ is equivalent to
the existence of a vector $w\in\mathbb{R}^n$ with positive coordinates and 
a number $\alpha\in[0,1)$ such that $Ew\leqslant\alpha w$.
More precisely, if such $w$ and $\alpha$ exist, 
then $r(E)\leqslant\alpha<1$ (see \cite[Corollary 8.1.29]{HJ13}).
Conversely, suppose that $r(E)<1$.
If $E$ is irreducible, 
then one can take $\alpha$ to be $r(E)$ and 
$w$ a Perron-Frobenius vector of $E$.
If $E$ is not irreducible, 
then consider the matrix $E'=(e_{ij}+\delta)$ for sufficiently small $\delta>0$,
and take $\alpha$ to be $r(E')$ and $w$ a Perron-Frobenius vector of $E'$.

Since $r(E^t)=r(E)$, we see that the condition $r(E)<1$ 
is also equivalent to the existence of a vector $w$ with positive coordinates 
and a number $\alpha\in[0,1)$ such that $E^t w\leqslant \alpha w$.

\begin{B}(\cite[Theorem 2.2]{F20})
(1) Suppose $w=(w_1,...,w_n)^t$ with positive coordinates and 
$\alpha\in[0,1)$ satisfy $Ew\leqslant\alpha w$.
Then for each $N\geqslant 1$,
$$
\|I-(I-A)^N-P\|\leqslant
(w_1+...+w_n)\max\{(1/w_1)\|P_1\|,...,(1/w_n)\|P_n\|\}\frac{\alpha^N}{1-\alpha}.
$$

(2) Suppose $w=(w_1,...,w_n)^t$ with positive coordinates and
$\alpha\in[0,1)$ satisfy $E^t w\leqslant \alpha w$.
Then for each $N\geqslant 1$,
$$
\|I-(I-A)^N-P\|\leqslant
(w_1\|P_1\|+...+w_n\|P_n\|)\max\{(1/w_1),...,(1/w_n)\}\frac{\alpha^N}{1-\alpha}.
$$
\end{B}

Using Theorem B, we can get concrete estimates for 
the rate of convergence of $I-(I-A)^N$ to $P$.
Suppose $E$ is irreducible and $r(E)<1$.
Take $\alpha$ to be $r(E)$ and $w$ a Perron-Frobenius vector of $E$.
Then we get
$$
\|I-(I-A)^N-P\|\leqslant
(w_1+...+w_n)\max\{(1/w_1)\|P_1\|,...,(1/w_n)\|P_n\|\}\frac{(r(E))^N}{1-r(E)}.
$$
Similarly, we can take $\alpha$ to be $r(E)$ and $w$ a Perron-Frobenius vector of $E^t$.
Then we get
$$
\|I-(I-A)^N-P\|\leqslant
(w_1\|P_1\|+...+w_n\|P_n\|)\max\{(1/w_1),...,(1/w_n)\}\frac{(r(E))^N}{1-r(E)}.
$$

\subsection{}
The assumption $r(E)<1$ is a \textit{sharp} sufficient condition 
for $X_1+...+X_n$ to be complemented in $X$.
More precisely, we have the following result.

\begin{C}(\cite[Theorem 2.3]{F20})
Let $E=(e_{ij})$ be an $n\times n$ matrix with 
$e_{ii}=0$ for $i=1,...,n$ and 
$e_{ij}\geqslant 0$ for any distinct $i,j\in\{1,...,n\}$.
If $r(E)=1$, 
then there exist a Banach space $X$, 
complemented subspaces $X_1,...,X_n$ of $X$ and 
projections $P_1,...,P_n$ onto $X_1,...,X_n$, respectively, such that
\begin{enumerate}
\item
$\|P_i x\|=e_{ij}\|x\|$, $x\in X_j$, for any distinct $i,j\in\{1,...,n\}$;
\item
$X_1,...,X_n$ are linearly independent;
\item
$X_1+...+X_n$ is closed and not complemented in $X$.
\end{enumerate}
\end{C}

\begin{remark}
If $r(E)>1$, the theorem can be applied to the matrix $(1/r(E))E$.
\end{remark}

\subsection{}
The aim of this paper is to obtain analogues of Theorems A, B, C
for the case when $X$ is a Hilbert space,
$X_1,...,X_n$ are closed subspaces of $X$,
and $P_1,...,P_n$ are orthogonal projections onto $X_1,...,X_n$, respectively.

\section{Results}

Let $X$ be a real or complex Hilbert space,
$X_1,...,X_n$ be closed subspaces of $X$,
and $P_1,...,P_n$ be orthogonal projections onto $X_1,...,X_n$, respectively.
As in the Banach space setting, we assume that nonnegative numbers 
$\varepsilon_{ij}$, $i\neq j$, $i,j\in\{1,...,n\}$ are such that
\begin{equation}\label{eq:H varepsilon_ij}
\|P_i x\|\leqslant\varepsilon_{ij}\|x\|,\quad x\in X_j,
\end{equation}
for any distinct $i,j\in\{1,...,n\}$.
Clearly, \eqref{eq:H varepsilon_ij} is equivalent to the 
inequality $\|P_i|_{X_j}\|\leqslant\varepsilon_{ij}$.

Now observe that for arbitrary closed subspaces $M$ and $N$ of 
the space $X$ $\|P_M|_N\|=\|P_N|_M\|$, 
where $P_M$ and $P_N$ are the orthogonal projections onto $M$ and $N$,
respectively.
Indeed, if $M=\{0\}$ or $N=\{0\}$, then $\|P_M|_N\|=\|P_N|_M\|=0$.
Assume that $M$ and $N$ are nonzero.
One can easily check that $(P_M|_N:N\to M)^*=P_N|_M:M\to N$.
Therefore $\|P_M|_N\|=\|P_N|_M\|$. 

Hence, we can and will assume that 
$\varepsilon_{ij}=\varepsilon_{ji}$ for any distinct $i,j\in\{1,...,n\}$.
Define the $n\times n$ matrix $E=(e_{ij})$ by
$$e_{ij}=
\begin{cases}
0, &\text{if $i=j$,}\\
\varepsilon_{ij}, &\text{if $i\neq j$.}
\end{cases}
$$
It is clear that $E$ is symmetric and nonnegative.
It follows that $r(E)$, the spectral radius of $E$, 
is the maximum eigenvalue of $E$.
Set $A:=P_1+...+P_n$.
Applying Theorem A and noting that
$$
\ker(P_1)\cap...\cap\ker(P_n)=
X_1^\bot\cap...\cap X_n^\bot=
(X_1+...+X_n)^\bot,
$$
we get the following result.

\begin{theorem}\label{Th:H}
If $r(E)<1$, 
then the subspaces $X_1,...,X_n$ are linearly independent and
their sum is closed in $X$.
Moreover, the sequence of operators $I-(I-A)^N$
converges uniformly to the orthogonal projection $P$ 
onto $X_1+...+X_n$ as $N\to\infty$.
\end{theorem}

\begin{remark}
A part on the closedness of Theorem~\ref{Th:H} follows from 
\cite[Theorem 1.2]{K11}.
\end{remark}

\begin{remark}
$r(E)<1$ $\Leftrightarrow$
the matrix $I-E$ is positive definite $\Leftrightarrow$
every leading principal minor of the matrix is positive.
For $n=2$ the inequality $r(E)<1$ is equivalent to
$\varepsilon_{12}^2<1$ $\Leftrightarrow$ $\varepsilon_{12}<1$.
For $n=3$ the inequality $r(E)<1$ is equivalent to
$$
\varepsilon_{12}^2+\varepsilon_{23}^2+\varepsilon_{31}^2+
2\varepsilon_{12}\varepsilon_{23}\varepsilon_{31}<1.
$$
\end{remark}

\begin{remark}
For the optimal choice $\varepsilon_{ij}=\|P_i|_{X_j}\|$, $i\neq j$,
the condition $r(E)<1$ can be formulated in terms of
the minimal angles between the subspaces.
Let us recall the definition of the minimal angle between two subspaces;
the notion was introduced by J.~Dixmier in \cite{D49}.
Let $M$ and $N$ be two closed subspaces of $X$.
Define the number $c_0(M,N)\in[0,1]$ by
$$
c_0(M,N)=\sup\{|\langle x,y\rangle|\,|\,
x\in M, \|x\|\leqslant 1, y\in N, \|y\|\leqslant 1\}.
$$
The minimal angle between $M$ and $N$ is the angle
$\varphi_0(M,N)\in[0,\pi/2]$ 
whose cosine is equal to $c_0(M,N)$.
It is well known, and one can easily prove, that
$c_0(M,N)=\|P_M|_N\|$, where $P_M$ is the orthogonal projection onto $M$.

Set $\varphi_{ij}:=\varphi_0(X_i,X_j)$, $i\neq j$, $i,j\in\{1,2,...,n\}$.
Then $\varepsilon_{ij}=\|P_i|_{X_j}\|=c_0(X_i,X_j)=\cos\varphi_{ij}$, 
$i\neq j$, and therefore
$$
e_{ij}=
\begin{cases}
0, &\text{if $i=j$,}\\
\cos\varphi_{ij}, &\text{if $i\neq j$.}
\end{cases}
$$
$r(E)<1$ $\Leftrightarrow$
the matrix $I-E$ is positive definite $\Leftrightarrow$
every leading principal minor of the matrix is positive.

For $n=2$ the inequality $r(E)<1$ is equivalent to
$\cos^2\varphi_{12}<1$ $\Leftrightarrow$ 
$\cos\varphi_{12}<1$ $\Leftrightarrow$
$\varphi_{12}>0$.
Thus, if $\varphi_{12}>0$, then $X_1\cap X_2=\{0\}$ and
the subspace $X_1+X_2$ is closed in $X$.
It is worth mentioning that the converse is also true, i.e.,
if $X_1\cap X_2=\{0\}$ and the subspace $X_1+X_2$ is closed in $X$,
then $\varphi_{12}>0$ (see \cite[Theorem 12]{D95}).

For $n=3$ the inequality $r(E)<1$ is equivalent to
$$
\cos^2\varphi_{12}+\cos^2\varphi_{23}+\cos^2\varphi_{31}+
2\cos\varphi_{12}\cos\varphi_{23}\cos\varphi_{31}<1.
$$
One can easily check that the last inequality is equivalent to
$\varphi_{12}+\varphi_{23}+\varphi_{31}>\pi$.
\end{remark}

To make the paper self-contained we will prove
Theorem~\ref{Th:H} in Section~\ref{S:proofs}.

For the rate of convergence of $I-(I-A)^N$ to $P$ 
we have the following estimate
which is more precise than that given by Theorem B.

\begin{theorem}\label{Th:H rate}
If $r(E)<1$, then $\|I-(I-A)^N-P\|\leqslant (r(E))^N$ 
for each $N\geqslant 1$.
\end{theorem}

The assumption $r(E)<1$ is a \textit{sharp} sufficient condition 
for $X_1+...+X_n$ to be closed in $X$.
More precisely, we have the following result.

\begin{theorem}\label{Th:H r(E)=1}
Let $E=(e_{ij})$ be a symmetric $n\times n$ matrix 
with $e_{ii}=0$ for $i=1,...,n$ and 
$e_{ij}\geqslant 0$ for any distinct $i,j\in\{1,...,n\}$.
If $r(E)=1$, 
then there exist a separable infinite dimensional Hilbert space $X$ and 
infinite dimensional closed subspaces $X_1,...,X_n$ of $X$ such that
\begin{enumerate}
\item
$\|P_i|_{X_j}\|=e_{ij}$ for any distinct $i,j\in\{1,...,n\}$; 
here $P_i$ is the orthogonal projection onto $X_i$, $i=1,...,n$.
\item
$X_1,...,X_n$ are linearly independent;
\item
$X_1+...+X_n$ is not closed in $X$.
\end{enumerate}
\end{theorem}

\begin{remark}
If $r(E)>1$, the theorem can be applied to the matrix $(1/r(E))E$.
\end{remark}

\begin{remark}
For a similar result see \cite[Theorem 1.2, part ``Moreover,...'']{K11}.
However, in \cite{K11} there is no proof for this part of Theorem 1.2.
\end{remark}

\section{Proofs}\label{S:proofs}

\subsection{Proof of Theorems \ref{Th:H} and \ref{Th:H rate}}

We first formulate the following lemma.

\begin{lemma}\label{L:H}
Let $H,K$ be Hilbert spaces and $T:H\to K$ be a continuous linear operator.
Suppose that
$$
\alpha\|x\|\leqslant\|Tx\|\leqslant\beta\|x\|,\quad x\in H,
$$
where $\alpha$ and $\beta$ are positive numbers.
Then
$$
\alpha\|y\|\leqslant\|T^* y\|\leqslant\beta\|y\|,\quad y\in Ran(T).
$$
\end{lemma}

The proof is simple and is omitted.

Now we are ready to prove Theorems \ref{Th:H} and \ref{Th:H rate}.
For simplicity of notation, we set $r:=r(E)$.
Let $X_1\oplus...\oplus X_n$ be 
the (orthogonal) direct sum of Hilbert spaces $X_1,...,X_n$.
Define the operator $S:X_1\oplus...\oplus X_n\to X$ by
$$
S(x_1,...,x_n)^t=x_1+...+x_n,\quad x_1\in X_1,...,x_n\in X_n.
$$
One can easily check that $S^*:X\to X_1\oplus...\oplus X_n$
acts as follows:
$$
S^*x=(P_1 x,...,P_n x)^t,\quad x\in X.
$$
Therefore $SS^*=P_1+...+P_n=A$.
For every $v=(v_1,...,v_n)^t\in X_1\oplus...\oplus X_n$ we have
\begin{align*}
&|\|Sv\|^2-\|v\|^2|=
|\|v_1+...+v_n\|^2-\|v_1\|^2-...-\|v_n\|^2|=\\
&=|\sum_{i\neq j}\langle v_i,v_j\rangle|\leqslant
\sum_{i\neq j}|\langle v_i,v_j\rangle|=
\sum_{i\neq j}|\langle P_i v_i,v_j\rangle|=\\
&=\sum_{i\neq j}|\langle v_i,P_i v_j\rangle|\leqslant
\sum_{i\neq j}\|v_i\|\|P_i v_j\|\leqslant
\sum_{i\neq j}\varepsilon_{ij}\|v_i\|\|v_j\|=\\
&=\langle E(\|v_1\|,...,\|v_n\|)^t,(\|v_1\|,...,\|v_n\|)^t\rangle\leqslant\\
&\leqslant r\|(\|v_1\|,...,\|v_n\|)^t\|^2=
r(\|v_1\|^2+...+\|v_n\|^2)=
r\|v\|^2.
\end{align*}
So $|\|Sv\|^2-\|v\|^2|\leqslant r\|v\|^2$, i.e.,
\begin{equation}\label{ineq:S}
(1-r)\|v\|^2\leqslant\|Sv\|^2\leqslant(1+r)\|v\|^2,\quad
v\in X_1\oplus...\oplus X_n.
\end{equation}

Since $\|Sv\|\geqslant\sqrt{1-r}\|v\|$, $v\in X_1\oplus...\oplus X_n$,
we conclude that $\ker(S)=\{0\}$ and $Ran(S)$ is closed in $X$
$\Rightarrow$
the subspaces $X_1,...,X_n$ are linearly independent and their sum
$X_1+...+X_n$ is closed in $X$.
Further, from \eqref{ineq:S} and Lemma~\ref{L:H} it follows that
$$
(1-r)\|x\|^2\leqslant\|S^* x\|^2\leqslant(1+r)\|x\|^2,\quad x\in Ran(S).
$$
We have
$$
\|S^* x\|^2=\langle S^*x,S^*x\rangle=\langle SS^*x,x\rangle=\langle Ax,x\rangle.
$$
The closed subspace $Ran(S)=X_1+...+X_n$ is invariant with respect to $A$.
Denote by $A'$ the restriction of $A$ to $X_1+...+X_n$.
Then we get
$$
(1-r)\|x\|^2\leqslant\langle A' x,x\rangle\leqslant(1+r)\|x\|^2,
\quad x\in X_1+...+X_n,
$$
and hence
$$
|\langle(A'-I)x,x\rangle|\leqslant r\|x\|^2,\quad x\in X_1+...+X_n.
$$
Since $A'-I$ is self-adjoint, we conclude that $\|A'-I\|\leqslant r$.

Let us estimate $\|I-(I-A)^N-P\|$, 
where $P$ is the orthogonal projection onto $X_1+...+X_n$.
Consider the orthogonal decomposition
$$
X=(X_1+...+X_n)\oplus(X_1+...+X_n)^\bot=
(X_1+...+X_n)\oplus (X_1^\bot\cap...\cap X_n^\bot).
$$
With respect to the decomposition we have
$P=I\oplus 0$ and $A=A'\oplus 0$.
Thus 
$$
I-(I-A)^N-P=-(I-A')^N\oplus 0
$$ 
and
$$
\|I-(I-A)^N-P\|=\|-(I-A')^N\|\leqslant\|A'-I\|^N\leqslant r^N\to 0
$$
as $N\to\infty$.

Theorems \ref{Th:H} and \ref{Th:H rate} are proved.

\subsection{Proof of Theorem \ref{Th:H r(E)=1}}

We will prove Theorem~\ref{Th:H r(E)=1} for the case when
the base field of scalars is $\mathbb{R}$;
the proof for $\mathbb{C}$ is similar.
For a number $\alpha\in(0,1)$ consider the matrix $I-\alpha E$.
Note that this matrix has the following properties:
\begin{enumerate}
\item
$I-\alpha E$ is a real symmetric matrix with diagonal elements equal to $1$;
\item
the least eigenvalue of this matrix is equal to $1-\alpha$.
Consequently, this matrix is positive definite.
\end{enumerate}
Therefore $I-\alpha E$ is the Gram matrix of some 
linearly independent collection of unit vectors of $\mathbb{R}^n$, say
$v^{(i)}=v^{(i)}(\alpha)$, $i=1,...,n$.
Let $L_i=L_i(\alpha)$ be the one-dimensional subspace 
spanned by $v^{(i)}$, $i=1,...,n$.
Denote by $P_i=P_i(\alpha)$ the orthogonal projection onto $L_i$, $i=1,...,n$.
Clearly, $P_i v=\langle v,v^{(i)}\rangle v^{(i)}$, $v\in\mathbb{R}^n$, $i=1,...,n$.
It follows that
$P_i v^{(j)}=\langle v^{(j)},v^{(i)}\rangle v^{(i)}=-\alpha e_{ij}v^{(i)}$
and consequently
$\|P_i|_{L_j}\|=\alpha e_{ij}$
for arbitrary $i\neq j$.
Further, since $1$ is an eigenvalue of $E$, 
there exists a unit vector $c=(c_1,...,c_n)\in\mathbb{R}^n$ such that $Ec=c$.
Then $(I-\alpha E)c=(1-\alpha)c$ and consequently 
$\langle (I-\alpha E)c,c\rangle=1-\alpha$.
We rewrite this equality as
$\sum_{i,j}\langle v^{(j)},v^{(i)}\rangle c_j c_i=1-\alpha,$
which is equivalent to
\begin{equation}\label{eq:norm of linear combination}
\|c_1 v^{(1)}+...+c_n v^{(n)}\|^2=1-\alpha.
\end{equation}

Now we are ready to construct a Hilbert space $X$ and 
its closed subspaces $X_1,...,X_n$ with the needed properties.
Take an arbitrary sequence $\alpha_k\in(0,1)$, $k=1,2,...$, 
which converges to $1$ as $k\to\infty$.
Set
$$
X=\bigoplus_{k=1}^\infty\mathbb{R}^n
$$
and
$$
X_i=\bigoplus_{k=1}^\infty L_i(\alpha_k),\quad i=1,...,n,
$$
where $\bigoplus$ is the (orthogonal) direct sum of Hilbert spaces.

First, let us show that $\|P_i|_{X_j}\|=e_{ij}$ for each pair $i\neq j$ 
(here $P_i$ is the orthogonal projection onto $X_i$).
It is clear that
$$
P_i=\bigoplus_{k=1}^\infty P_i(\alpha_k).
$$
Hence,
$$
\|P_i|_{X_j}\|=
\sup\{\|P_i(\alpha_k)|_{L_j(\alpha_k)}\|\,|\,k=1,2,...\}=
\sup\{\alpha_k e_{ij}\,|\,k=1,2,...\}=
e_{ij}.
$$

Further, since the vectors 
$v^{(1)}(\alpha_k),...,v^{(n)}(\alpha_k)$ are linearly independent
for $k=1,2,...$, we conclude that the subspaces $X_1,...,X_n$
are linearly independent.

It remains to show that $X_1+...+X_n$ is not closed in $X$.
Suppose that $X_1+...+X_n$ is closed in $X$.
Let $X_1\oplus...\oplus X_n$ be 
the (orthogonal) direct sum of Hilbert spaces $X_1,...,X_n$.
Define the operator $S:X_1\oplus...\oplus X_n\to X$ by
$$
S(x_1,...,x_n)^t=x_1+...+x_n,\quad x_1\in X_1,...,x_n\in X_n.
$$
Since $X_1,...,X_n$ are linearly independent, we conclude that
$\ker(S)=\{0\}$.
Moreover, $Ran(S)=X_1+...+X_n$ is closed in $X$.
These properties imply that $S$ is an isomorphic embedding.
Hence, there exists a number $\beta>0$ such that
$\|Su\|\geqslant\beta\|u\|, u\in X_1\oplus...\oplus X_n.$
We rewrite this inequality as
\begin{equation}\label{ineq:norm of sum}
\|x_1+...+x_n\|^2\geqslant
\beta^2(\|x_1\|^2+...+\|x_n\|^2),\quad x_1\in X_1,...,x_n\in X_n.
\end{equation}
Now we choose
$x_i=(0,...,0,c_i v^{(i)}(\alpha_k),0,0,...)^t, i=1,...,n.$
By~\eqref{eq:norm of linear combination}
and~\eqref{ineq:norm of sum} 
we get $1-\alpha_k\geqslant\beta^2$.
But $\alpha_k\to 1$ as $k\to\infty$ and thus we get a contradiction.
Hence, $X_1+...+X_n$ is not closed in $X$.

Theorem~\ref{Th:H r(E)=1} is proved.

\subsection*{Acknowledgements}

The research was funded by 
Institute of Mathematics of NAS of Ukraine.
This research was supported by 
the Project 2017-3M from the Department of Targeted Training of
Taras Shevchenko National University of Kyiv at the NAS of Ukraine.

\end{document}